\documentclass[preprint,12pt,authoryear]{elsarticle}

\usepackage{amssymb}

\usepackage{graphicx,psfrag,epsf}

\journal{Journal of Statistical Planning and Inference}
\newcommand{\beit}{\begin{itemize}}
\newcommand{\eit}{\end{itemize}}

\newcommand{\bi}{\begin{itemize}}
\newcommand{\ei}{\end{itemize}}

\newcommand{\vekk}[1]{}

\newcommand{\bes}{\begin{frame}}
\newcommand{\ens}{\end{frame}}

\newcommand{\bx}{{\mathbf{x}}}
\newcommand{\bQ}{{\mathbf{Q}}}

\newcommand{\notes}[1]{}
\newcommand{\delete}[1]{}
\newcommand{\gtcomment}[1]{}
\newcommand{\st}{|}
\newcommand{\ceps}{{\cal E}} 
\renewcommand{\Pr}{\mbox{Pr}} 
\newcommand{\into}{\mbox{$\: \rightarrow \:$}} 
\begin{document}

\begin{frontmatter}



\title{On the proper treatment of improper distributions\tnoteref{label1}}

\tnotetext[label1]{Accepted manuscript published as: Lindqvist B. H., Taraldsen G., On the proper treatment of improper
distributions. \textit{J. Statist. Plann. Inference} (2017), https://doi.org/10.1016/j.jspi.2017.09.008}


\author{Bo H. Lindqvist}
\ead{bo.lindqvist@ntnu.no}

\author{Gunnar Taraldsen}
\ead{gunnar.taraldsen@ntnu.no}
\address{Department of Mathematical Sciences, Norwegian University of Science and Technology, N-7491 Trondheim, Norway}

\begin{abstract}
The axiomatic foundation of probability theory presented by Kolmogorov has been the basis of modern theory for probability and statistics. In certain applications it is, however, necessary or convenient to allow improper (unbounded) distributions, which is often done without a theoretical foundation. The paper reviews a recent theory which includes improper distributions, and which is related to Renyi's theory of conditional probability spaces.  It is in particular demonstrated how the theory leads to simple explanations of apparent paradoxes known from the Bayesian literature. Several examples from statistical practice with improper distributions are discussed in light of the given theoretical results, which also include a recent theory of convergence of proper distributions to improper ones.
\end{abstract}


\end{frontmatter}


\section{Introduction}
\label{intro}

Bayes' formula forms the basis of Bayesian statistics. Suppose a parameter $\theta$ is of interest, and that we have data $x$ which is supposed to give information about $\theta$. The idea of Bayesian inference is to first express one's prior knowledge (some would call it \textit{uncertainty}) of $\theta$ in the form of a  \textit{prior distribution}, commonly given in the form of a density function $\pi (\theta)$, and then combine this knowledge with the new knowledge provided by the data~$x$. The influence of  $\theta$ on the data $x$ is modeled by a  \textit{statistical model}, represented by the conditional density of the data given the parameter, $f (x \st \theta)$. Note that $f(x|\theta)$ will sometimes be interpreted as the 'likelihood' of $\theta$ for a given observation $x$, in which case the function $\theta \mapsto f(x|\theta)$ will be the \textit{likelihood function}. 

Bayes' formula is  used to express the updated information about $\theta$ obtained after $x$ is observed, given in the form of the  \textit{posterior distribution,}
\begin{equation}
\label{eqBay1}
\pi (\theta \st x) = \frac{f(x \st \theta) \pi(\theta)}{f (x)}.
\end{equation}
Here $f (x) = \int f(x \st \theta) \pi (\theta) d\theta$
is the \textit{marginal density} of $x$.
Equation~(\ref{eqBay1}) is one version of Bayes' theorem. 

The algorithm for calculation of posterior distributions given by (\ref{eqBay1})
is well defined as long as 
$0 < f (x) < \infty$, and in this case it will always lead to a proper probability distribution in the sense that $\pi(\theta|x)$ is a non-negative  function which integrates to~1. 

The usual proof of Bayes' formula
is restricted to the case where
$\pi (\theta)$ is a probability density, where basic rules of probability are used in the derivation. As pointed out above we may however get a proper distribution as the output of the formula even if the prior $\pi(\theta)$ is not proper. This fact is the obvious excuse for using Bayes' formula for improper priors. 

A natural question is now, why should one want to use improper distributions for $\theta$? In practice, improper distributions often result from the search for so-called \textit{non-informative} priors. The most prominent example of such priors is the Jeffreys' prior, which is proportional to the square root of the determinant of the Fisher information matrix and has the key property of being invariant under reparameterizations. 

Intuitively, a non-informative prior should be one that does not favor any parameter values above others, suggesting ``flat'' priors. In practice, this often means to use proper standard probability models such as normal, gamma or uniform distributions with (very) large variances. Taking the limit as the variance tends to $\infty$ is then most often the excuse for using an improper prior density which is constant over the complete parameter space. Such distributions may, however, not have the invariance properties required by Jeffreys' priors. It is in fact well-known that flat priors may be very informative on non-location parameters. We refer to \cite{bernardo} for an interesting discussion on non-informative and improper priors. 

As already indicated, posterior distributions computed by Bayes' formula are proper probability distributions only under the condition of $0 < f(x)<\infty$. 
Standard Bayesian calculations are, however, made by just recognizing the proportionality
\[
\pi (\theta \st x) \propto  f(x \st \theta) \pi(\theta),
\]
which can be used without loss of information in the case that Bayes' formula gives a proper distribution, but gives an improper posterior distribution in case $f(x)$ is not finite. The latter case, if ignored, may lead to misleading inferences, as discussed later in this paper.

Improper distributions also appear naturally in certain non-Bayesian analyses.  \cite{lindtar}
considered conditional sampling of data $x$ given a sufficient statistic $T(x)$ for the unknown parameter $\theta$, which has numerous applications in statistical inference \citep{casella}. The key  is that this conditional distribution is independent of the value of $\theta$. The  general idea of the conditional sampling method of \cite{lindtar} was to use this fact, but instead of fixing the value of $\theta$, to let it be a random quantity with some suitable distribution. Under certain mild restrictions, this distribution can be freely chosen, often with improper distributions giving the most efficient methods. Improper distributions appear likewise as useful ingredients in fiducial statistics, see for example \cite{tarlind-ann} and \cite{taraldsen-fiducial}.

The purpose of the present paper is to review and discuss some important aspects of the use of improper priors in  statistical practice. Some would say that no new theory is needed, since improper priors are just approximations to proper ones. As will be discussed in the paper, this is a too simple attitude. The literature on Bayesian statistics includes a lot of paradoxes and misleading conclusions due to improper priors and posteriors. There are, however, not a lot of theoretical treatments of proper versus improper distributions in the literature. Some exceptions are, e.g., \cite{hartigan}, \cite{pollard} and the more recent paper by \cite{mccullagh}.   

Our point of departure will be the paper by \cite{TAS2010} which has a slightly different view than the above references. The idea is here simply to allow infinite probabilities in Kolmogorov's axioms. While this implies that all random variables have infinite mass, all conditional distributions will be proper probability distributions under a certain crucial condition which turns out to be equivalent to the above mentioned condition of finite $f(x)$. Formally, this condition is the $\sigma$-finiteness of the random quantity that is conditioned on, here $x$. Details will be given in Section 2 which reviews the theoretical results of \cite{TAS2010}.
 
The above idea is not new, however. We quote from \cite{Renyi1962}, motivating the introduction of improper distributions: 
\begin{quote}
``\textit{One can indeed give an axiomatic theory of probability which matches
the above-mentioned requirements.
This theory contains the theory of Kolmogorov as a special case. The
fundamental concept of the theory is that of conditional probability; it
contains cases where ordinary probabilities are not defined at all.}''\end{quote}

In a footnote he adds: \begin{quote}``\textit{The idea of such a theory is due to Kolmogorov himself;
he, however, did not publish anything about it.}''\end{quote}

The theory presented in \cite{TAS2010} is in fact closely related to Renyi's theory of conditional probability spaces \citep{renyi1970}. This connection is studied in more detail in \cite{taraldsen-renyi}.

Having introduced the basic elements of the theory of \cite{TAS2010} in Section~\ref{theory}, we proceed to Section~\ref{bayesian} which discusses some consequences of the theory when applied to Bayesian statistics. In particular we investigate in some detail a so-called marginalization paradox presented by \cite{stone1972}.  Section~\ref{gibbs} is devoted to Gibbs sampling, where a possible pitfall is the fact that posteriors may be improper even if all full conditionals are proper. A recent theoretical paper on approximation of improper priors, \cite{bioche}, is briefly reviewed in Section~\ref{convergence}. This is an important paper giving precise conditions  for convergence of proper priors to improper ones and for convergence of the corresponding posterior distributions. Section~\ref{intrinsic} discusses a class of improper models which is popular in spatial statistics. Some concluding remarks are finally given in Section~\ref{concluding}.

\section{The theoretical framework}
\label{theory}
\subsection{The modified Kolmogorov axioms}
As in Kolmogorov's axioms we consider an abstract space $\Omega$ of outcomes, where events $A$ are represented by subsets of $\Omega$ and where the family $\ceps$ of events is assumed to be a $\sigma$-algebra.
We next let the measurable space $(\Omega,\ceps)$ be equipped with a
fixed law~$\Pr$ with
\beit
\item
$\Pr(A) \ge 0$ 
for all $A \in \ceps$. \item 
$\Pr (A_1 \cup A_2 \cup \cdots) = \Pr (A_1) + \Pr (A_2) + \cdots$ \\
whenever $A_1,A_2,\ldots$ are pairwise disjoint events.
\eit
 However, where Kolmogorov adds the axiom $\Pr (\Omega) = 1$, we  assume only
 \beit
 \item $\Pr(\mbox{\O})=0$,
 \eit 
 and hence allow the case  $\Pr (\Omega) = \infty$. Note that the above axioms are exactly the axioms of a positive measure from standard measure theory \citep{royden}.

\subsection{Random quantities}
\label{sec2.2}
A random quantity $X$ with values in a measurable space $(\Omega_X,\ceps_X)$, 
is identified with a measurable function
$X: \Omega \into \Omega_X$, i.e., such that  $(X \in A) \equiv \{\omega \in \Omega \; \st \; X (\omega) \in A\}$ is an event in $\ceps$ 
 for any event $A$ in $\ceps_X$.
The law $\Pr$ on $\Omega$ now induces the law $\Pr_X$ of a random quantity $X$
 by defining
\[
\Pr_X (A)  = \Pr(X \in A) \mbox{ for  $A \in \ceps_X$}.
\]
Hence the joint law of a pair $(X,Y)$  is determined by 
\[
\Pr_{X,Y}(A \times B) = \Pr((X,Y) \in A \times B) \mbox{ for $A \in \ceps_X$, $B \in \ceps_Y$},
\]
while marginal laws are found from
\[
\Pr_X (A) =  \Pr_{X,Y} (A \times \Omega_Y) \mbox{ for $A \in \ceps_X$}. 
\]

The random quantity  $Y$ is called $\sigma$-finite if
the law $\Pr_Y$ is $\sigma$-finite, i.e., if
there exist events $E_1, E_2, \ldots \in \ceps_Y$  with:
\[ \Omega_Y = \cup_i E_i \mbox{ and } \Pr_Y (E_i) < \infty \mbox{ for $i = 1,2,\ldots$} \]

\subsection{Conditional distributions}
\label{sec2.3}
 
A key feature of our approach is that if $Y$ is $\sigma$-finite, then we can define a \textit{unique proper conditional probability }
\[ \Pr^y (A) \equiv \Pr(A | Y=y) \] as a function of $y \mbox{ for each } A \in \ceps$. The following approach equals the standard approach for definition of conditional probabilities and expectation in ordinary probability theory. 

For a given event $A$ in $\ceps$, conditional probabilities  should satisfy, for all $B \in \ceps_Y$:
 \begin{eqnarray}
 \Pr(A \cap (Y \in B))
 &=& \int_B \Pr(A |Y=y)\Pr_{Y} (dy)  \nonumber \\ &=& \label{royden}
\int_B \Pr^y (A) \,\Pr_{Y} (dy) . 
\end{eqnarray}
By the assumed $\sigma$-finiteness of the measure $\Pr_Y$, the  Radon-Nikodym theorem  \citep{royden}
states exactly that the function 
$g(y) = \Pr^y (A)$ exists and is uniquely (a.e.) defined by the above.

Since the measure $\Pr_Y$  must satisfy $\Pr_Y(Y \in B) = \int_B \Pr_Y(dy)$ for all $B \in \ceps_Y$, it is seen by letting $A=\Omega$ in (\ref{royden}) and using uniqueness of $g(y)$, that we have
\[
\Pr^y (\Omega) = 1.
\]
Under regularity conditions, which we will not pursue here, we may from this conclude that conditional laws $\Pr^y$ can always be represented as proper probability distributions, as long as $Y$ is $\sigma$-finite. If, on the other hand, $Y$ is 
not $\sigma$-finite, then $\Pr^y$ is not defined due to the requirement of $\sigma$-finiteness in Radon-Nikodym's theorem.

Having defined the conditional law $\Pr^y$ on $\Omega$, we now define the conditional distribution of a random quantity $X$ given $Y=y$ for $A \in \ceps_X$ by
\[
   \Pr_X^y (A) = \Pr^y(X \in A).
   \]

\subsection{A Bayesian statistical model}
A Bayesian statistical model involves
an observation, represented by a random quantity $X :\Omega \into \Omega_X$, 
and a random parameter $\theta$, represented as a $\sigma$-finite random quantity 
$\Theta :\Omega \into \Omega_\Theta$. The law $\Pr_\Theta$ 
of $\Theta$ is then the prior distribution.

The conditional distribution of $X$ given $\Theta=\theta$, i.e.,
$\{\Pr_X^\theta : \theta \in \Omega_\Theta \}$, defines in a consistent way a \textit{statistical model}. This follows directly from the above approach since $\Theta$ is assumed to be $\sigma$-finite and since conditional distributions are always proper.

\subsection{Implications of improper prior $\Pr_\Theta$}
So far we have not specified the value of $\Pr(\Omega)$. Suppose that
$\Theta $ is $\sigma$-finite with $\Pr_\Theta(\Omega_\Theta)=\infty$. We claim that this implies that $\Pr$ is $\sigma$-finite and $\Pr(\Omega)=\infty$. To see this, suppose $A_1, A_2, \ldots \in \ceps_\Theta$  are such that
\[ \Omega_\Theta = \cup_i A_i \mbox{ and } \Pr_\Theta (A_i) < \infty \mbox{ for $i = 1,2,\ldots$} \]
Then
\[
\Omega = (\Theta \in \cup_i A_i)=\cup_i (\Theta \in A_i),
\]
which implies that $\Pr$ is also by necessity 
improper and $\sigma$-finite. This follows since $\Pr(\Omega) = \Pr_\Theta(\Omega_\Theta) = \infty$ and 
$\Pr(\Theta \in  A_i) = \Pr_\Theta (A_i) < \infty$.

Assume now that $\Pr(\Omega)=\infty$. Then every random quantity has an improper law, since 
\[
\Pr_X(\Omega_X)= \Pr\{\omega:X(\omega) \in \Omega_X\} = \Pr(\Omega)=\infty.
\]

On the other hand, a random quantity $X$ is not necessarily  $\sigma$-finite, even if $\Pr$ is. 
Namely, let $X$ take values $0$ and $1$. Then
\[
\infty  = \Pr(X \in \Omega_X) = \Pr(X=0) + \Pr(X=1) 
\] 
and at least one of these is necessarily equal to $\infty$. Hence $X$ is not $\sigma$-finite.

\subsection{Bayesian posteriors}
Recall that  a Bayesian model is given by a $\sigma$-finite law $\Pr_\Theta$, the prior distribution, and an observation $X$  with distribution $\Pr_X^\theta$. For an observation $x$, Bayesian inference considers the \textit{posterior law}, i.e., the conditional law of $\Theta$ given $X=x$, which in our notation is
$
 \Pr_\Theta^x 
 $.
This conditional distribution is well defined if  $X$ is $\sigma$-finite, in which case it is a \textit{proper} probability distribution. On the other hand, if  $X$ is \textit{not} $\sigma$-finite, then the posterior is not defined.
Hence, in the current theory there is nothing such as an \textit{improper posterior}!

\texttt{\texttt{}}

\section{Bayesian statistics and marginalization paradoxes}
\label{bayesian}
\subsection{The absolutely continuous case}
Random quantities are said to be absolutely continuous if they can be defined by densities with respect to Lebesgue measure, for example  $\Pr_{X,Y} \sim  f(x,y)$. 
The \textit{marginal density} of $X$ is then given by the density 
$f (x)  = \int f(x,y) dy $, where $f(x)=\infty$ is a permitted value.
 
 It is  seen that $X$ with density $f(x)$ is $\sigma$-finite according to the definition of Section~\ref{sec2.2} if and only if $f(x) < \infty$ (a.e.), and the approach of the previous section can be shown to lead to the
Bayes' formula (\ref{eqBay1}), 
which corresponds to $\Pr_\Theta^x$ in the notation of  Section~\ref{sec2.3}. 

\subsection{What may ``go wrong'' with improper distributions?} 
A  prior model for $\theta=(\theta_1,\theta_2)$  in Bayesian statistics is commonly given on the form of a joint density $\pi(\theta_1,\theta_2)=\pi_1(\theta_1)\pi_2(\theta_2)$ of the pair $(\Theta_1,\Theta_2)$, where $\pi_1(\theta_1)$, $\pi_2(\theta_2)$ are two non-negative, finite-valued functions. One would then say that the parameters are given ``independent priors, with marginal priors  $\pi_1(\theta_1)$ and $\pi_2(\theta_2)$''. In practice, one might have chosen one or both of the ``marginal'' priors $\pi_1(\theta_1)$ and $\pi_2(\theta_2)$ as improper ones. To be concrete, suppose $\pi_1(\theta_1)$ is a proper probability density, while $\pi_2(\theta_2)$ is improper, i.e., integrates to $\infty$. As indicated above, it would be tempting to call $\pi_1(\theta_1)$ and  $\pi_2(\theta_2)$ the \textit{marginal} densities of $\Theta_1$ and $\Theta_2$, respectively. But are they?  

By the definition given in Section~\ref{sec2.2}, the marginal density of $\Theta_1$ is $\int \pi(\theta_1,\theta_2)d\theta_2= \pi_1(\theta_1)\int \pi_2(\theta_2)d\theta_2$, which however equals $\pi_1(\theta_1) \cdot \infty$ since $\pi_2(\cdot)$ is improper. Since this equals $\infty$ whenever $\pi_1(\theta_1)>0$, it follows that $\pi_1(\cdot)$ is \textit{not} the marginal density of $\Theta_1$! However, integrating instead with respect to $\theta_1$ and recalling that $\pi_1(\cdot)$ was assumed to be a probability density, we find that the marginal density of $\Theta_2$ is $ \int \pi(\theta_1,\theta_2)d\theta_1= \pi_2(\theta_2)\int \pi_1(\theta_1)d\theta_1 = \pi_2(\theta_2)$, showing that $\pi_2(\cdot)$ is indeed the marginal density of~$\Theta_2$.  

So which interpretation can we give of $\pi_1(\theta_1)$? Since the marginal density of $\Theta_2$, $\pi_2(\theta_2)$, is finite (although not proper), we conclude that $\Theta_2$ is $\sigma$-finite. Hence it has meaning to condition on it, and it can be seen that $\pi_1(\theta_1)$ is the \textit{conditional density} of $\Theta_1$ given $\Theta_2=\theta_2$, using the approach of Section~\ref{sec2.3}. In particular this shows that even if $\Theta_1$ is not $\sigma$-finite, and has an infinite marginal density, it has a well defined proper conditional distribution given the $\sigma$-finite random quantity $\Theta_2$.

\subsection{A marginalization paradox \citep{stone1972}}
\label{sec33}
For given parameters $\Theta=\theta,\Phi=\phi$, let $X$ and $Y$ be independent and exponentially distributed with hazard rates, respectively, $\theta \phi$ and $\phi$.  Suppose the interest is in the ratio $\theta$ of the hazard rates, which suggests consideration of the ratio $Z=Y/X$.

Let the joint prior distribution of $(\Theta,\Phi)$ be given by
 $\pi(\theta,\phi) =\pi(\theta) \cdot 1$, where $\pi(\theta)$ is  proper. (Note that by the previous subsection, $\pi(\theta)$ is not the \textit{marginal} density of $\theta$.)
The joint density of $(X,Z,\Theta,\Phi)$ is readily obtained to be
\begin{equation}
\label{init}
f(x,z,\theta,\phi) = \theta \phi^2 x e^{-\phi x(\theta+z)}\pi(\theta).
\end{equation}

Integration with respect to $(\theta,\phi)$ shows that $(X,Z)$ is $\sigma$-finite, and we readily get the marginal conditional distribution of $\Theta$ given $(X,Z)$ to be
\begin{equation}
\label{cross}
f(\theta|x,z) \propto \frac{\theta \pi(\theta)}{(\theta+z)^3}.
\end{equation}
Since this does not depend on $x$, it is tempting to conclude that the right hand side of (\ref{cross}) is also the conditional distribution of $\Theta$ given $Z=z$, i.e.,
\begin{equation}
\label{fz1}
f(\theta|z) = f(\theta|x,z) \propto \frac{\theta \pi(\theta)}{(\theta+z)^3}.
\end{equation}
   
   Starting differently, by integrating out $x$ in (\ref{init}) and conditioning with respect to $(\theta,\phi)$ (which is obviously $\sigma$-finite), we get
\[
f(z|\theta,\phi) = \int_0^\infty \theta \phi^2 x e^{-\phi x(\theta+z)}dx = \frac{\theta}{(\theta+z)^2}.
\]
Since this depends only on $\theta$, one might suggest that 
$
f(z|\theta) = \frac{\theta}{(\theta+z)^2}.
$
and from this obtain
\begin{equation}
\label{fx2}
f(\theta|z) \propto f(z|\theta)\pi(\theta) = \frac{\theta \pi(\theta)}{(\theta+z)^2}.
\end{equation}
But (\ref{fz1}) and (\ref{fx2}) contradict each other! It is therefore not clear how to proceed if one wants to do inference on $\theta$ based on $Z$ alone. This is an example of a \textit{marginalization paradox.} 

So what is the problem? Considering the approach of Section~\ref{theory}, the problem is that the marginal distribution of $Z$ is not $\sigma$-finite, so one is not allowed to condition on it. The conclusion in (\ref{fz1}) is hence not correct. In fact, neither is the distribution of $\Theta$ $\sigma$-finite, making the conclusion (\ref{fx2}) incorrect as well.

The clue is that we have above, in fact twice, used the generally invalid result that
\[
f(\theta|x,z) \mbox{ does not depend on $x$ } \Rightarrow f(\theta|z)=f(\theta|x,z).
\]
This is well-known to hold for probability distributions, but holds for improper distributions only provided $(X,Z)$ and $Z$ both have $\sigma$-finite distributions.

\subsection{Marginalization paradox revisited}
\label{34}
In order to understand better the mechanisms of the previous example, let us redefine the problem and let the prior of $(\theta,\phi)$ be given by
\[
   \pi(\theta,\phi)= \pi(\theta)h(\phi),
\]
where $h(\phi)$ may be proper or improper, while $\pi(\theta)$ is proper as before, unless otherwise stated below.
Multiplying (\ref{init}) by $h(\phi)$ and integrating with respect to $\phi$ we get
\begin{equation}
\label{xzt}
  f(x,z,\theta)= \frac{\theta \pi(\theta)}{x^2} \int_0^\infty
  u^2 e^{-u(\theta+z)}h\left(\frac{u}{x}\right)du.
\end{equation}
The joint marginal distribution of $(z,\theta)$ is obtained by integrating with respect to $x$, which gives
\begin{eqnarray}
   f(z,\theta) &=& \theta \pi(\theta) \int_0^\infty u^2e^{-u(\theta+z)}
   \left[ \int_0^\infty \frac{1}{x^2} h\left(\frac{u}{x}\right)dx \right] du \nonumber \\
  &=& \frac{\theta \pi(\theta) }{(\theta+z)^2} \int_0^\infty h(w)dw. \label{nils}
\end{eqnarray}
Hence, if $h$ is proper, then for inference about $\theta$ we may base ourselves on $Z$ only and use the relation (\ref{fx2}). This case corresponds to the typical frequentist approach for this example, where one concludes that the distribution of $Z$ depends on the parameters only via $\theta$. 
 In the case where $h$ is improper, we get however $f(z,\theta) = \infty$ in (\ref{nils}), and neither $Z$ nor $\Theta$ are $\sigma$-finite. 
 
 Recalling that we want to make inference about $\theta$, let us go back to (\ref{xzt}). It follows that 
 \begin{equation}
\label{doublestar}
  f(\theta|x,z) \propto \frac{\theta \pi(\theta)}{(\theta +z)^3}
  \int_0^\infty w^2e^{-w}h\left(\frac{w}{x(\theta+z)} \right)dw.
\end{equation}
 Setting $h \equiv 1$ leads to the right hand side of (\ref{fz1}). Strange enough, choosing $h$ to be the improper density $h(\phi)=1/\phi$, it follows that
 \begin{equation}
 \label{DD}
 f(\theta|x,z) \propto \frac{\theta \pi(\theta)}{(\theta+z)^2}
\end{equation}
 in which case we would apparently not have a marginalization paradox. Still, however,  $Z$ is not $\sigma$-finite, so we cannot conclude that (\ref{DD}) equals $f(\theta|z)$. 
It is notable that \cite{stone1972} explain this apparent absense of a marginalization paradox by the fact that we
now use the prior $1/\phi$ for $\Phi$ which is the common prior for a scale parameter.  In our 
  opinion this seems to be more like a coincidence since we here use an improper~$h$ under which $Z$ is  not $\sigma$-finite.
	
	As another comment on (\ref{DD}), note that the pair $(X,Z)$ is  $\sigma$-finite even if we let $\pi(\theta)$ be the improper density $\pi(\theta)=1/\theta$, while we keep $h(\phi)=1/\phi$. This is seen by integrating (\ref{xzt}) with respect to $\theta$. Hence (\ref{DD}) is meaningful and leads after normalization to the posterior density for $\theta$ given by
	\begin{equation}
	\label{freq}
	f(\theta|x,z) = \frac{z}{(\theta + z)^2}.
	\end{equation}
	It is interesting to note that the density (\ref{freq}) also appears as the optimal invariant confidence distribution for $\theta$ in a frequentist approach involving the observations $X$ and $Y$, where $\theta$ is the parameter of interest. The argument follows Schweder and Hjort (2016), Chapter 5.

 We close the present section by returning to the original assumption where  $h(\phi) \equiv 1$ and $\pi(\theta)$ is unspecified, but proper.  Consider a proper density which can be seen as an approximation to the improper $h(\phi)$, e.g.,
 \begin{equation}
 \label{hM}
    h_M(\phi) = \frac{1}{M} I(0<\phi \leq M),
 \end{equation}
 where $I(\cdot)$ is the indicator function and $M>0$ is considered to be large. From (\ref{doublestar}) we get
 \begin{eqnarray}
      f(\theta|x,z) & \propto&  \frac{\theta \pi(\theta)}{(\theta+z)^3}\int_0^{xM(\theta+z)} w^2e^{-w}dw \nonumber \\
 &=&  \frac{\theta \pi(\theta)}{(\theta+z)^3}\left[ 2 - e^{-A} (A^2+2A+2) \right], \label{xM}
 \end{eqnarray}
 where  $A = xM(\theta+z)$. 
 
\begin{figure}[h!]
\centering
\includegraphics[width=12cm,angle=0]{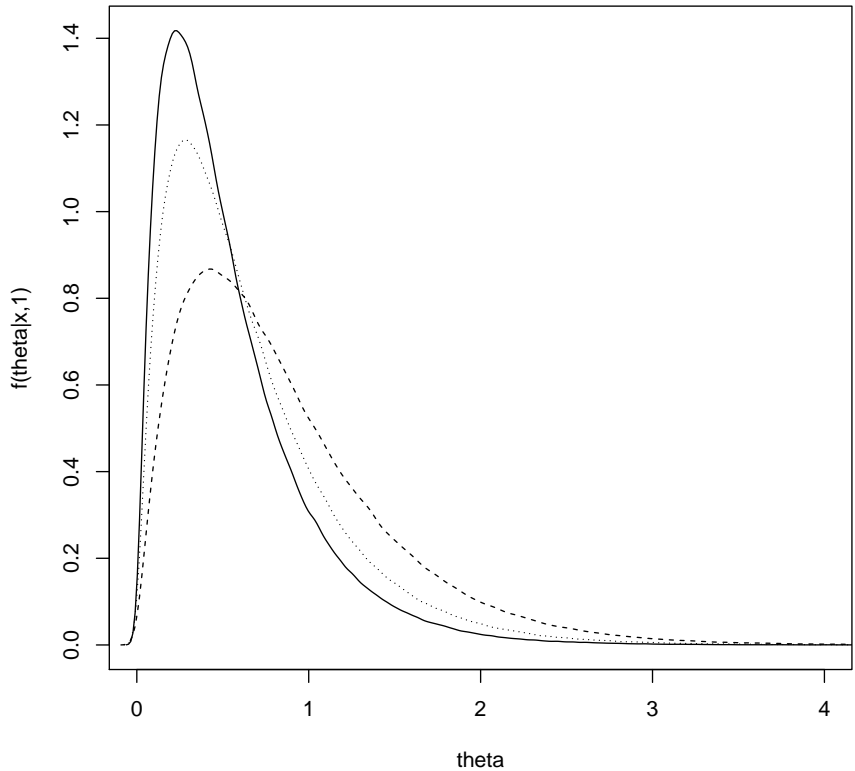}  
\caption{Marginalization paradox example. Let $\pi(\theta)=\exp(-\theta)$ and suppose $z=1$ is observed. Solid line: the (normalized) density (\ref{xM}) with $x=1, M=500$ (which is indistinguishable from (\ref{cross})).  Dashed line: the (normalized) density (\ref{xM}) with $x=0.001, M=500$. Dotted line: the (normalized) density (\ref{DD}).}
\label{3posterior}
\end{figure}

 It is seen that the limit as $M$ tends to infinity in (\ref{xM}) is consistent with (\ref{cross}). This is in fact a consequence of a general result in \cite{bioche} (see Section~\ref{convergence}), since $h_M \rightarrow h$ in their approach with $h \equiv 1$. As seen from (\ref{xM}), the convergence of $f(\theta|x,z)$ as $M$ tends to infinity is not uniform in the observations $(x,z)$. This point was made by \cite{akaike} in his discussion of certain marginalization paradoxes. More precisely, Akaike questioned the common interpretation of an improper  prior distribution as a limit of  proper prior distributions, and he argued that
 an improper prior  can more adequately  be described as the limit of certain data adaptive proper prior
 distributions. He concluded that a prior distribution without data adaptability may produce poor inference due
 to a gross misspecification of the prior. We illustrate this point in Figure~\ref{3posterior}. It is seen that, even if we set $M$ as a large number (here 500) in (\ref{hM}),  the posterior distribution for $\Theta$ in (\ref{xM}) depends rather distinctly on the value of $x$, at least  for  small $x$. This is despite the non-appearance of $x$ in the right hand side of (\ref{cross}).  For higher values of $x$, like $x=1$, the posterior is however indistinguisable from the one we get when $M$ equals infinity.  The figure also includes a corresponding plot of (\ref{fx2}) which illustrates the difference between (\ref{cross}) and (\ref{fx2}).

\vekk{ 
 Further remarks,
 
 $\pi(\theta)$ is not the marginal distribution of $\theta$, since marginal is $\infty$.
 
 Note: A frequentist might be tempted to use only Z, since it has a distribution that depends only on $\theta$. Question: What about Lehmann-Scheffe, etc. Do we have an optimal inference there? $(X,Z)$ is minimal sufficient (and complete), but 
 
 The distribution of $\theta$ given $x,z$ depends on both x and z when $h(\phi)$ is proper,
 \[
   f(\theta|x,z) \propto \frac{\theta \pi(\theta)}{(\theta +z)^2}
   \int_0^\infty w^2e^{-w}h\left(\frac{w}{x(\theta+z)} \right)dw
   \] 
 
 This is the best expression! Can be used also when $h \equiv 1$.
 
 $Z/\theta$ is Fisher(2,2)
 } 
 
\section{An example from Gibbs-sampling}
\label{gibbs}

\subsection{Gibbs-sampling from improper posterior distribution}

\cite{hobert} gave an example showing that the output from Gibbs sampling corresponding to an improper posterior distribution may still appear perfectly reasonable. The authors' advice is thus that before implementing a Gibbs chain one should check that the posterior is proper. For this it is important to note that  propriety of the conditionals of a Gibbs chain does not imply that the full posterior is proper (see example below). 

\cite{gelfandsahu1999} consider similar problems with Gibbs sampling, focusing on parameter identifiability and posterior propriety. In particular, they provide rather general 
results for propriety of posteriors in the case of GLMs. 
As a simple illustration they consider in an earlier technical report \citep{gelfandsahu1996} the following example. 

Let $Y|\theta_1,\theta_2 \sim N(\theta_1+\theta_2,1)$ with improper prior $\pi(\theta_1,\theta_2)=1$.  
Then the joint distribution of $(Y,\Theta_1,\Theta_2)$ is
\[
f(y,\theta_1,\theta_2)=
(1/\sqrt{2\pi})e^{-(1/2)(y-\theta_1-\theta_2)^2},
\]
  leading to the marginal density of $Y$ given as
\[ 
\int \int f(y,\theta_1,\theta_2) d\theta_1 \; d\theta_2 = \infty.
\]
Hence $Y$ is not $\sigma$-finite, so the posterior $f(\theta_1,\theta_2|y)$ does not exist (or, is \textit{improper}).
  
On the other hand, the pairs $(Y,\Theta_1)$ and  $(Y,\Theta_2)$ are both $\sigma$-finite, so 
the following conditional distributions exist and are proper:
\begin{eqnarray}
\Theta_1|\theta_2,y & \sim & N(y-\theta_2,1) \label{t1},\\
\Theta_2|\theta_1,y & \sim & N(y-\theta_1,1) \label{t2}.
\end{eqnarray}

\begin{figure}[h!]
\centering
\includegraphics[width=6.5cm,angle=0]{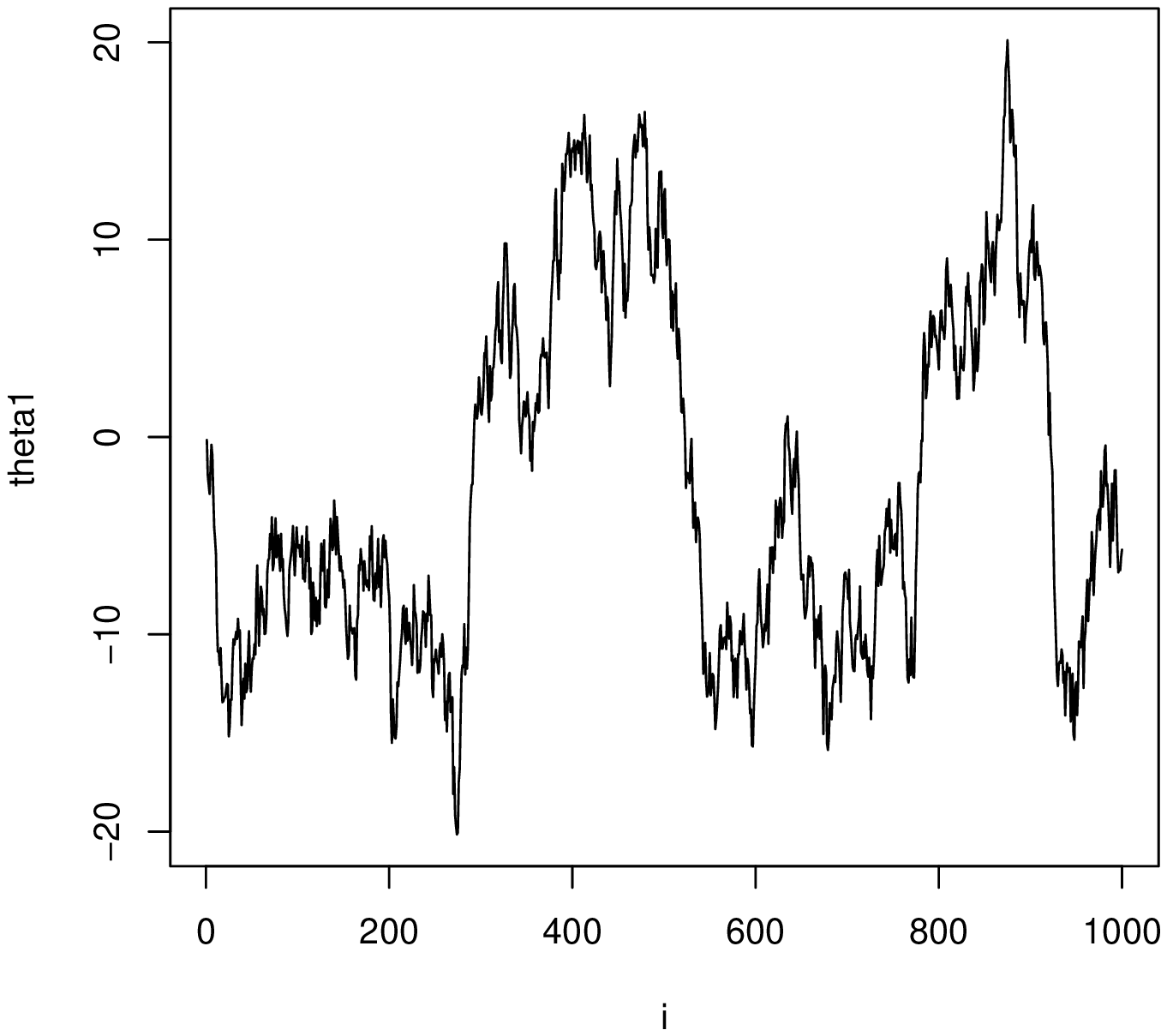}, \includegraphics[width=6.5cm,angle=0]{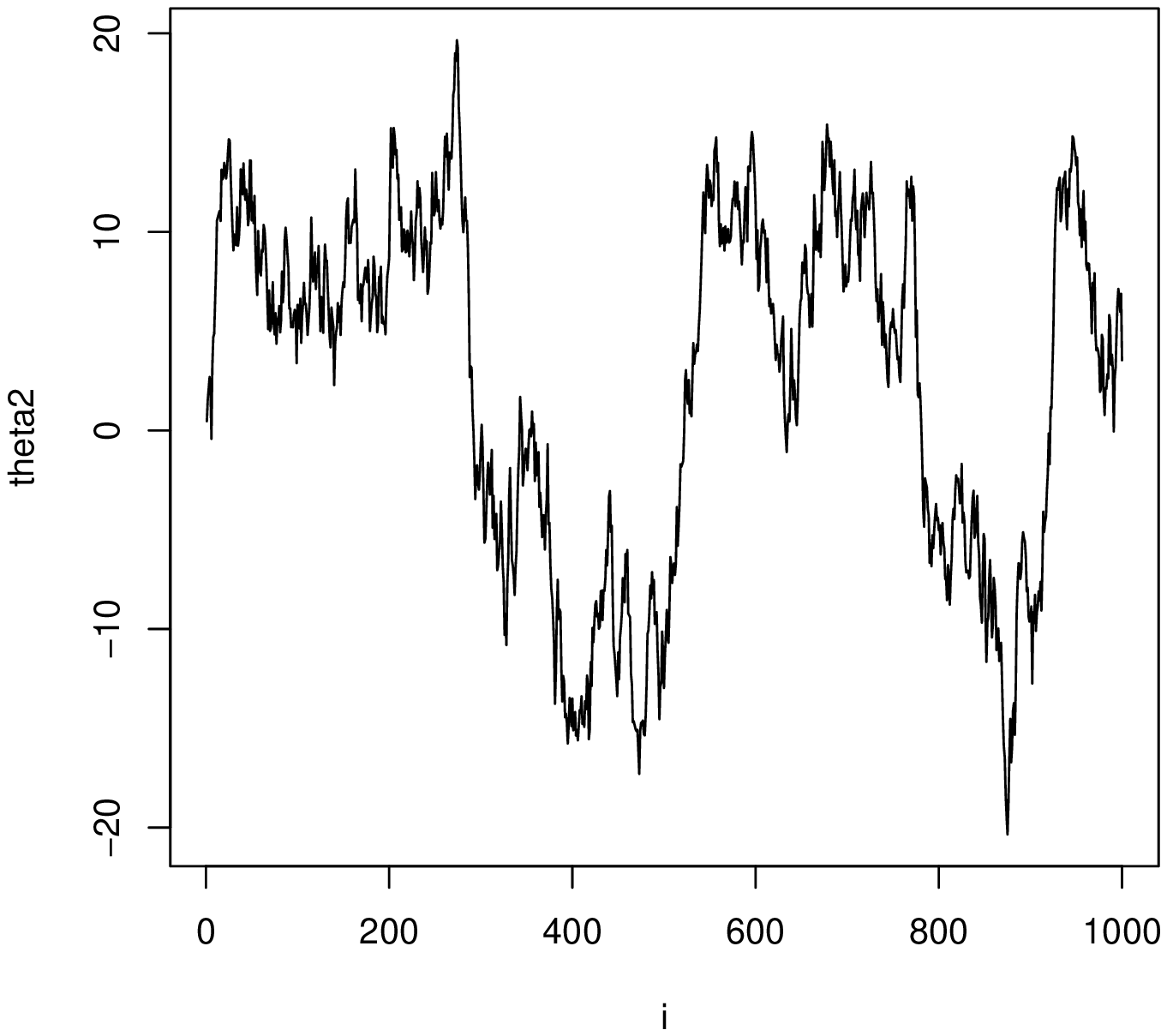}
\caption[]{Gibbs chains for $\theta_1$ (left) and $\theta_2$ (right), drawn from (\ref{t1}) and (\ref{t2}), respectively.}
\label{Gibbsfig1}
\end{figure}

Thus Gibbs-sampling of pairs $(\theta_1,\theta_2)$ for given $y$ is possible.
The question is, however, how the pairs $(\theta_1,\theta_2)$ will behave, knowing that the posterior $f(\theta_1,\theta_2|y)$ does not exist. Figure~\ref{Gibbsfig1} shows a simulation from (\ref{t1})-(\ref{t2}). The large fluctuations seen in the plots are due to the impropriety of the joint posterior given $y$. 

  \begin{figure}[h!]
\centering
\includegraphics[width=6.5cm,angle=0]{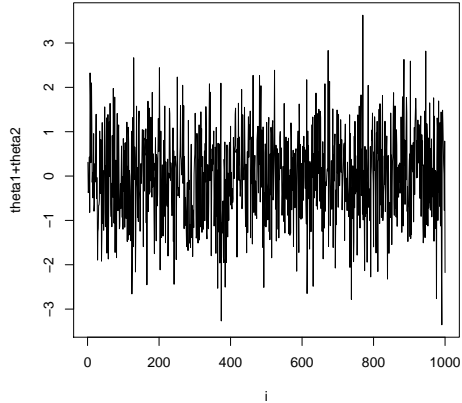}
\caption[]{Simulated values of $\delta = \theta_1 + \theta_2$ using (\ref{t1}) and (\ref{t2}).}
\label{Gibbsfig2}
\end{figure}

\subsection{The proper embedded posterior \citep{gelfandsahu1999}}
Gelfand and Sahu observed, however, that if one makes a 1-1 transformation 
 \[ (\theta_1,\theta_2) \into (\delta,\rho), \; \;  \mbox{ where $\delta=\theta_1+\theta_2$,}\]  then the distribution $\delta|y \sim N(y,1)$ can be recovered in the Gibbs-sampling. Indeed, the plot of $\delta=\theta_1+\theta_2$ from (\ref{t1}) and (\ref{t2}) (Figure~\ref{Gibbsfig2}), is apparently well-behaved. Gelfand and Sahu (1999) call $\delta|y$  the unique proper embedded  posterior, regarding it as embedded within the improper posterior for $(\theta_1,\theta_2)$. 

A critical remark is of course appropriate here in view of the previous theory. Since $Y$ is not $\sigma$-finite it has apparently no meaning to consider $\delta | y$. On the other hand, we clearly have from (\ref{t2}) that $(\Theta_1+\Theta_2|y,\theta_1) \sim N(y,1)$, i.e., if we let $\rho \equiv \theta_1$, we have $(\delta|y,\rho) \sim N(y,1)$.  Thus Gelfand and Sahu's conclusion is similar to the one that we deemed to be incorrect in connection with the marginalization paradox. Namely, when $Y$ is not $\sigma$-finite, then even if the density of $(\delta|y,\rho)$ does not depend on $\rho$, this is not the conditional density of $\delta$ given $y$.

 The  nice behavior of $\delta$ in the simulation can be explained theoretically as follows. Suppose that the prior distribution of $(\Theta_1,\Theta_2)$ is given by $\pi(\theta_1,\theta_2)=g(\theta_1)$, where $g(\theta_1)$ is 
 a proper density. 
  Then under the transformation $(\theta_1,\theta_2) \rightarrow (\rho,\delta)$, we have
 \begin{equation}
  \label{tytti}
f(y,\rho,\delta)=
(1/\sqrt{2\pi})e^{-(1/2)(y-\delta)^2} \cdot g(\rho).
\end{equation}
 In this model $Y$ is clearly $\sigma$-finite (in fact, the marginal density of $Y$ is the constant 1). 
 Thus the posterior $ \pi(\rho,\delta|y)$ exists and is given by (\ref{tytti}). The marginal posterior of $\delta$ is hence
 \[
 \delta|y \sim N(y,1)
 \]
 whatever be the density $g(\rho)$, as long as it is proper. \cite{gelfandsahu1999} let $g(\rho)$ correspond to $N(0,\tau^2)$ and let $\tau^2 \rightarrow \infty$, and concluded that also in the limit will have $\delta|y \sim N(y,1)$. This is, however, not a valid conclusion since $Y$ is now not $\sigma$-finite.

\subsection{Using  proper priors for both $\theta_1$ and $\theta_2$}
A proper posterior for $(\theta_1,\theta_2)$ can of course be achieved by giving $(\Theta_1,\Theta_2)$ a proper prior. Assume for example  that $\Theta_1$ and $\Theta_2$ are independent with
$\Theta_1 \sim N(0,\tau^2)$, $\Theta_2 \sim N(0,\kappa^2)$. Then, as shown in Gelfand and Sahu (1996),
\begin{eqnarray*}
   \Theta_1 | \theta_2,y &\sim& N\left( \frac{\tau^2}{1+\tau^2} (y-\theta_2),\frac{\tau^2}{1+\tau^2} \right) ,\\
   \Theta_2 | \theta_1,y &\sim& N\left( \frac{\kappa^2}{1+\kappa^2} (y-\theta_1),\frac{\kappa^2}{1+\kappa^2} \right) .
   \vekk{ 
   \theta_1 | ,y &\sim& N\left( \frac{\tau^2 y}{1+\tau^2+\kappa^2},\frac{\tau^2(1+\kappa^2)}{1+\tau^2+\kappa^2} \right)   \\
    \theta_2 | ,y &\sim& N\left( \frac{\kappa^2 y}{1+\tau^2+\kappa^2},\frac{\kappa^2(1+\tau^2)}{1+\tau^2+\kappa^2} \right) .  } 
   \end{eqnarray*}    
    If $\tau$ and $\kappa$ are large, then the trajectories of
the Gibbs chains for $\theta_1$  and $\theta_2$, respectively,  will still tend to drift in a way similar to the behavior in Figure~\ref{Gibbsfig1}. Thus, if we use a proper but diffuse priors for $\Theta_1$ and $\Theta_2$, the posteriors will be proper but will in practice be
indistinguishable from those obtained under the corresponding limiting improper prior. As concluded by \cite{gelfandsahu1996}, an implicit byproduct of this observation is the infeasibility of numerical sampling based
diagnostics for propriety of posteriors. A similar conclusion is expressed by \cite{hobert}.

\section{Convergence of priors and posteriors (Bioche and Druilhet)}
\label{convergence}

\subsection{$q$-vague convergence of measures}

\cite{bioche} propose a convergence mode for measures allowing  a sequence
of probability measures to have an improper limiting measure. They also study  convergence of corresponding posterior distributions. 
 
Technically the authors  study the set of positive Radon-measures on the state space $\Omega_\Theta$, i.e., the set of positive measures $\Pi$ which are finite on compact subsets of $\Omega_\Theta$. Noting that the output of Bayes' formula (\ref{eqBay1}) is unchanged if $\pi(\theta)$ is multiplied by a constant, they define the equivalence relation $\Pi \sim \Pi'$ to mean that there is an $\alpha >0$ such that $\Pi=\alpha \Pi'$. Their basic space of measures is then the corresponding quotient space, equipped with the quotient topology resulting from vague convergence of positive Radon measures. Convergence in this topology has been denoted as q-vague convergence. A similar quotient topology is introduced by \cite{taraldsen-renyi}.   

A useful way of expressing the definition of $q$-vague convergence is the following: A sequence of positive Radon-measures $\{\Pi_n\}$ converges $q$-vaguely to $\Pi$ if there exists a sequence $\{a_n\}$ such that 
\[ 
a_n\Pi_n \rightarrow \Pi \mbox{ (vaguely),}
\]  
   (see, e.g., \cite{billingsley} for the definition of vague convergence).
	
From this definition it is not difficult to prove that 
for any improper distribution $\Pi$ there is a sequence $\Pi_n$ of proper distributions such that 
$\Pi_n \rightarrow \Pi$ ($q$-vaguely). In this case, the $a_n$ given in the above definition  tend to $\infty$ as $n$ increases. 

As an example, consider the proper distribution with density $h_M$ given by (\ref{hM}). We claim that $h_M \rightarrow h$ (q-vaguely) as $M \rightarrow \infty$, where $h \equiv 1$. To see this we need to find constants $a_M$ such that $a_Mh_M \rightarrow h$ (vaguely) as $M \rightarrow \infty$, i.e., such that
\[
\int a_M h_M(\phi) f(\phi) d\phi \rightarrow \int f(\phi)d\phi
\]
for each continuous function with compact support. But this is clear by the dominated convergence theorem by letting $a_M = M$ for all $M>0$.  

\cite{bioche} also consider convergence of posterior densities. If $f(x|\theta)$ is the likelihood of the data $x$ and $\pi(\theta)$ is the prior density, then they define the posterior distribution as the distribution with density in the equivalence class corresponding to $ \pi(\theta|x) \propto  f(x|\theta)\pi(\theta)$,
thus allowing also improper posterior distributions. 

Their main proposition on convergence of posteriors states that if 
for the priors we have  $\pi_n \rightarrow \pi$ ($q$-vaguely), and if $\theta  \mapsto f(x|\theta)$ is continuous, then the posteriors converge in the sense that $\pi_n(\cdot|x) \rightarrow \pi(\cdot|x) $ ($q$-vaguely). We have already seen an example in Section~\ref{34}, where the posterior distribution  (\ref{xM}) converges  to (\ref{cross}) as $M$ tends to infinity. Note, however, that while the question of uniform convergence in $x$ was made a point in our example, this issue is not considered by \cite{bioche}.

At first glance it seems that the above cited result on convergence of posteriors justifies the common excuse for using improper priors, namely that they are limits of proper priors and hence that the posteriors are limits of posteriors based on proper priors. However, we have already seen problems connected to such a view. Next we shall see another type of misinterpretation of improper limits of proper distributions, which in turn may give completely misleading results regarding posterior distributions.

\subsection{The Jeffreys-Lindley paradox (Bioche and Druilhet)}
Let $X|\theta \sim N(\theta,1)$, and consider testing of the null hypothesis $H_0:\theta=0 \mbox{ versus the alternative } H_1: \theta \ne 0$. 
Suppose we have a prior distribution for $\theta$ given by
\[
\pi(\theta) = \frac{1}{2} \delta_0 + \frac{1}{2} I(\theta \ne 0),
\]
 where $\delta_0$ is a point mass at $\theta=0$ and $I(\cdot)$ is the indicator function. This means that we have a prior belief of 1/2 in $H_0$, while the remaining probability 1/2 is distributed according to Lebesgue measure on $H_1$. A straightforward calculation gives  
\[
    \pi(0|x) = \left( 1+\sqrt{2 \pi}e^{x^2/2} \right)^{-1}
\]
implying $\pi(0|x) \le \left( 1+\sqrt{2 \pi} \right)^{-1} \approx 0.285$ whatever be the data $x$. 
    
  Using instead the proper prior measure 
  \[
  \pi_n(\theta) = \frac{1}{2} \delta_0 + \frac{1}{2} N(0,n^2)
  \] 
  we get
  \[
    \pi_n(0|x) = \left( 1+\sqrt{\frac{1}{1+n^2}}e^{\frac{n^2x^2}{2(1+n^2)}} \right)^{-1} .
    \]
But this converges to 1 as $n \rightarrow \infty$, in conflict with the above calculation which was based on an apparently equivalent argument using the limiting prior. The result has therefore been considered as a paradox. 

The clue, as presented by \cite{bioche}, is that while $N(0,n^2)$ converges q-vaguely to Lebesgue measure on the real line, the measure $\frac{1}{2} \delta_0 + \frac{1}{2} N(0,n^2)$ converges to $\frac{1}{2} \delta_0 \sim  \delta_0$ and not to $\frac{1}{2} \delta_0 + 
\mbox{Lebesgue}$ measure which one might believe.  This explains the paradox, noting that by the convergence result for posteriors,  the limiting posterior is a point mass at 0 as well.

\section{Intrinsic Gaussian Markov random fields (IGMRF)}
\label{intrinsic}

Intrinsic conditional autoregressions (ICAR) are widely used in  spatial statistics and dynamic
linear models \citep{besag1991}. These models are improper versions of conditional autoregressive models (CAR) as introduced as spatial models by \cite{besag1974}. Important special cases of CAR and ICAR models are  Gaussian Markov Random Fields (GMRF) and the intrinsic (improper) versions denoted IGMRF, see \cite{rue} for a thorough treatment including applications.

As discussed by \cite{lavine}, the fact that the intrinsic models correspond to improper distributions, implies that care should be taken in their use and interpretation. 

\subsection{The first order random walk}
Following \cite{rue} we use this simple special case of an IGMRF to illustrate some of the main issues regarding IGMRF models. 
 
Let $\bx= (x_1,x_2,\ldots,x_n)$ be the successive observations of a random walk, assuming independent increments
\[
   \Delta x_i = x_{i+1} - x_{i} \sim_{iid} N(0,\kappa^{-1}), \; i=1,2,\ldots,n-1.
   \]
   The IGMRF model specifies the density of $\bx$  to be the density obtained from these increments (only), giving
     \begin{eqnarray}
      f(\bx|\kappa) &\propto& \kappa^{(n-1)/2} \exp\left(-\frac{\kappa}{2} \sum_{i=1}^{n-1} (\Delta x_i)^2\right) \nonumber \\
    &=&  \kappa^{(n-1)/2} \exp\left(-\frac{\kappa}{2} \bx^T\bQ \bx \right). \label{fxk}
      \end{eqnarray}
Here  the \textit{structure matrix} $\bQ$  (displayed in \cite{rue}, p. 96) is positive semi-definite, with exactly one eigenvalue equal to 0, implying that $f(\bx|\kappa)$ is an improper density. 
   
 Statistical inference in models involving IGMRFs may involve making inference about the precision parameter $\kappa$. In a Bayesian analysis, one must  typically assign to $\kappa$ a  hyperprior and work with (\ref{fxk}) as a likelihood function. In this connection, \cite{lavine} question the use of the constant $\kappa^{(n-1)/2}$ appearing in (\ref{fxk}). In the following discussion let us replace (\ref{fxk}) by 
  \begin{equation}
 \label{star}
      f(\bx|\kappa) \propto  c(\kappa) \exp\left(-\frac{\kappa}{2} \bx^T\bQ \bx \right), 
      \end{equation} 
    thus making the appropriate choice of $c(\kappa)$ the main issue. As reported by \cite{lavine}, this choice has been discussed in several papers during the last two decades. \cite{besag1991} in fact used $c(\kappa)= \kappa^{n/2}$, which was used by WinBUGS (\cite{winbugs}) until it was changed to $c(\kappa)= \kappa^{(n-1)/2}$ following derivations appearing in, e.g., \cite{knorr2003} and  \cite{hodges2003}. 
		
   Rue and Held (2005) justify the density (\ref{fxk}) as follows: Consider first the 1-1 transformation  $(x_1,x_2,\ldots,x_n) \leftrightarrow (\Delta x_1,\ldots,\Delta x_{n-1},\bar x)$ where $\bar x$ is the average of the $x_i$. Assuming that $\bx$ is multivariate normal, $(\Delta x_1,\ldots,\Delta x_{n-1})$ and $\bar x$ are stochastically independent. We may hence write down the following proper density for $\bx$, indexed by $k$ for the purpose of later taking limits,
  \begin{equation}
  \label{hrue}
   \tilde f_k(\bx|\kappa) = f(\bx|\kappa) \cdot \breve f_k(\bar x)
  \end{equation}
  Here $f(\bx|\kappa)$ is the density (\ref{fxk}) while $\breve f_k(\bar x)$ is normal with zero expectation and precision $\gamma_k > 0$. Suppose $\gamma_k \rightarrow 0$ as $k \rightarrow \infty$.  We may then invoke Proposition 2.15 of \cite{bioche} to show that 
  \begin{equation}
  \label{kkk}
   \tilde f_k(\bx|\kappa) \rightarrow f(\bx|\kappa) \mbox{ as } k \rightarrow \infty
     \end{equation}  
  The interpretation of this is that the improper density (\ref{fxk}) is the limit of a sequence of proper densities for $\bx$. This derivation can also be interpreted as adding to the model for the $\Delta x_i$ a prior specification for $\bar x$ given in the form of a constant prior. 
	
	\cite{lavine} point, however, to a problem with this conclusion, having to do with the non-uniqueness of marginals in cases involving improper distributions and related to our discussion in Section~\ref{theory}. 
	\vekk{ 
	By the definition of $f(\bx|\kappa)$ in  (\ref{fxk}), and by the representation of $\bx$ as $(\Delta x_1,\ldots,\Delta x_{n-1},\bar x)$ one might suggest that $f(\bx|\kappa)$ is the marginal distribution of $(\Delta x_1,\ldots,\Delta x_{n-1})$. However, replacing $\breve f_k(\bar x)$ in (\ref{hrue}) by 1 and integrating   over $\bar x$ gives $\infty$, which is hence the formal marginal distribution of $\bx$ as discussed in Section cxz.
	}
To illustrate, essentially following Lavine and Hodges, we consider the modified 1-1 transformation of $\bx$ given as
\[
   (\Delta x_1,\ldots,\Delta x_{n-1}, \bar x \Delta x_1).
\]
  It follows by the ordinary transformation formula (involving a Jacobi-determinant), starting from (\ref{hrue}), that we have
 \[
   \tilde f_k(\bx|\kappa) = f(\bx|\kappa) \cdot \breve f_k(\bar x/\Delta x_1)\cdot \frac{1}{|\Delta x_1|}.
  \]
 Again, letting $\gamma_k \rightarrow 0$, we
 get  
 \[
   \tilde f_k(\bx|\kappa) \rightarrow f(\bx|\kappa)\cdot \frac{1}{|\Delta x_1|} \mbox{ as } k \rightarrow \infty,
     \]
thus giving a limit different from (\ref{kkk}).

\cite{lavine}  conclude that essentially all the arguments given in the literature for the value of the constant $c(\kappa)$ in some way are flawed. Their conclusion is therefore that any value of this constant may do. This is of course also in accordance with the previous section where the quotient topology for distributions was used, and where improper (as well as proper) distributions were identified with equivalence classes only. 

Having said this, there seem to be good reasons to use the form (\ref{fxk}). It follows from \cite{rue}, p. 90-91, who considered a more general case, that (\ref{fxk}) when restricted to $\bx$ such that $\bar x = \mu$, is the conditional density of $\bx$ given $\bar x = \mu$. Here $\mu$ can be any real number, but it seems that $\mu=0$ is commonly used. Furthermore, the specification of $\mu$ enables one to simulate from the distribution (\ref{fxk}) (see \cite{rue}, p. 92).

\subsection{Bayesian analysis with IGMRFs}

In a Bayesian inference with $\kappa$ as a parameter we  consider $f(\bx|\kappa)$ in (\ref{star}) as a likelihood function. It should then be noted that $f(\bx|\kappa)$ is improper and hence not proportional to a proper distribution, which is the case for commonly considered likelihood functions.

Let $\pi(\kappa)$ be the prior density of $\kappa$, possibly improper. 
The natural definition of the joint distribution of $(\bx,\kappa)$ is then
\begin{equation}
\label{twostar}
    f(\bx,\kappa) = f(\bx|\kappa) \pi(\kappa).
    \end{equation}
Thus the marginal density of $\kappa$ is
   \[
    \int f(\bx,\kappa)d\bx  = \int f(\bx|\kappa) \pi(\kappa)d\bx = \infty,
    \] 
so  $\pi(\kappa)$ is in fact not the marginal distribution of $\kappa$. 
But still, by the theory of Section~\ref{theory}, the posterior density $\pi(\kappa|\bx)$ is well defined provided $\bx$ is $\sigma$-finite. This holds if the integral over $\kappa$ of (\ref{twostar}) is finite for 
(almost) all $\bx$, i.e., if
\[
       \int \pi(\kappa) c(\kappa) \exp\left(-\frac{\kappa}{2} \bx^T\bQ \bx \right) d\kappa < \infty.
\]
A sufficient condition for this is clearly that $\int \pi(\kappa) c(\kappa)  d\kappa < \infty$. The conclusion of the above is that Bayesian inference for $\kappa$ is well-behaved under reasonable restrictions, as soon as the constant $c(\kappa)$ has been determined.

\section{Concluding remarks}
\label{concluding}

In this paper we have presented, and discussed in view of several examples, a simple theoretical approach which enables the inclusion of improper priors in Bayesian analyses. A special feature of the approach is that both parameters and observations
are represented as random quantities defined on a common underlying space $\Omega$. The clue has been to allow the probability $\Pr$ in Kolmogorov's axioms to be a $\sigma$-finite law with $\Pr (\Omega) = \infty$. In fact it was shown in Section~\ref{theory} that $\Pr (\Omega) = \infty$ is \textit{necessary} if improper priors are to be included.
    
What makes this a sensible theory is the fact that all conditional distributions, given $\sigma$-finite random quantities, are proper distributions. In particular this property leads to a consistent treatment of statistical models and a theoretically based condition for posterior propriety. 

The relation to Renyi's theory of conditional probability spaces has been mentioned earlier. In this connection we would also like to quote from \cite{lindley}. In the Preface to his classical test on probabilities he writes:
\begin{quote}
\textit{ The axiomatic structure used here is not the usual one associated with the
name of Kolmogorov. Instead one based on the ideas of Renyi has been
used. The essential dfference between the two approaches is that Renyi's
is stated in terms of conditional probabilities, whereas Kolmogorov's is in
terms of absolute probabilities, and conditional probabilities are defined
in terms of them. Our treatment always refers to the probability of
A, given B, and not simply to the probability of A. In my experience
students benefit from having to think of probability as a function of two
arguments, A and B, right from the beginning. The conditioning event,
B, is then not easily forgotten and misunderstandings are avoided. These
ideas are particularly important in Bayesian inference where one's views
are influenced by the changes in the conditioning event.}
\end{quote}
    
\section*{References}









\end{document}